\documentclass[12pt]{article} 
\usepackage{amsfonts,amsmath,latexsym,amssymb,mathrsfs,amsthm}
\usepackage{slashbox}
\usepackage{caption}
\usepackage[colorlinks=true]{hyperref}

\let\OLDthebibliography\thebibliography
\renewcommand\thebibliography[1]{
  \OLDthebibliography{#1}
  \setlength{\parskip}{3pt}
  \setlength{\itemsep}{0pt plus 0.3ex}
}

\def\numberlikeadb{\global\def\theequation{\thesection.\arabic{equation}}}
\numberlikeadb
\newtheorem{theorem}{Theorem}[section]

\newtheorem{corollary}[theorem]{Corollary}

\newtheorem{remark}[theorem]{Remark}

\usepackage{lscape}
\usepackage{caption}
\usepackage{multirow}
\begin{document}

\title{Inequalities for the modified Bessel function of the second kind and the kernel of the Kr\"{a}tzel integral transformation}
\author{Robert E. Gaunt\footnote{School of Mathematics, The University of Manchester, Manchester M13 9PL, UK}
}

\date{January 2017} 
\maketitle

\vspace{-10mm}

\begin{abstract}We obtain new inequalities for the modified Bessel function of the second kind $K_\nu$ in terms of the gamma function.  These bounds follow as special cases of inequalities that we derive for the kernel of the Kr\"{a}tzel integral transformation.  
\end{abstract}


\noindent{{{\bf{AMS 2010 Subject Classification:}}} Primary 33C10  

\section{Introduction}

The modified Bessel function of the second kind $K_\nu$ is an important and widely used special function.  There exists a substantial literature concerning inequalities for the modified Bessel function of the second; see, for example, \cite{luke} and \cite{baricz2} and references therein.  In a recent work, \cite{gaunt} derived the following simple lower bound for the function $K_0$:
\begin{equation}\label{koo}\frac{1}{\sqrt{x+\frac{1}{2}}}<\frac{\Gamma(x+\frac{1}{2})}{\Gamma(x+1)}<\sqrt{\frac{2}{\pi}}\mathrm{e}^xK_0(x).
\end{equation} 
%The approach used by \cite{gaunt} involved exploiting a useful integral representation of the modified Bessel function $K_0$, which, more generally, for the function $K_\nu$ is given by (\cite{olver}, formula 10.32.8):
In this note, we generalise this inequality to the modified Bessel function $K_\nu$ for $\nu\geq0$.  In deriving our inequality, we follow the approach of \cite{gaunt} by exploiting the following integral representation of the modified Bessel function of the second kind (\cite{olver}, formula 10.32.8):
\begin{equation}\label{krep}K_\nu(x)=\frac{\sqrt{\pi}(\frac{1}{2}x)^{\nu}}{\Gamma(\nu+\frac{1}{2})}\int_1^{\infty}\mathrm{e}^{-xt}(t^2-1)^{\nu-\frac{1}{2}}\,\mathrm{d}t, \quad x>0.
\end{equation}  
This integral representation of the modified Bessel function $K_\nu$ closely resembles the kernel
\begin{equation}\label{lrep}\lambda_\nu^{(n)}(x)=\frac{(2\pi)^{(n-1)/2}\sqrt{n}\big(\frac{x}{n}\big)^{n\nu}}{\Gamma(\nu+1-\frac{1}{n})}\int_1^\infty(t^n-1)^{\nu-\frac{1}{n}}\mathrm{e}^{-xt}\,\mathrm{d}t, \quad \text{$\nu>\frac{1}{n}-1$, $n=1,2,\ldots$,} 
\end{equation}
of the Kr\"{a}tzel integral transformation \cite{kratzel} (see also \cite{bar} and references therein for further properties) defined by
\begin{equation*}\mathcal{L}_\nu^{(n)}\{f\}(z)=\int_0^\infty\lambda_\nu^{(n)}(zt)f(t)\,\mathrm{d}t, \quad \mathrm{Re}\,x>0.
\end{equation*}
Indeed, a simple manipulation yields the relation
\begin{equation}\label{rel}\lambda_\nu^{(2)}(x)=2\bigg(\frac{x}{2}\bigg)^\nu K_\nu(x).
\end{equation}

Due to the similarity between the representations (\ref{krep}) and (\ref{lrep}), our approach to bounding the kernel $\lambda_\nu^{(n)}$ is no more difficult than bounding the modified Bessel function $K_\nu$.  In this note, we exploit the representation (\ref{lrep}) to derive inequalities for the kernel $\lambda_\nu^{(n)}$ and then use (\ref{rel}) to immediately deduce inequalities for the modified Bessel function $K_\nu$.  These inequalities are a natural generalisation of the inequality (\ref{koo}).

\section{Results and proofs}

The following is the main result of this note.

\begin{theorem}\label{beso} (i). Let $x>0$.  Then, for $0\leq \nu\leq\frac{1}{n}$, $n=1,2,\ldots$, we have
\begin{equation}\label{label}\lambda_{\nu}^{(n)}(x)\geq (2\pi)^{(n-1)/2}\frac{\sqrt{n}}{n-1}\bigg(\frac{x}{n}\bigg)^{n\nu}\frac{\Gamma(\frac{x}{n-1}+\frac{1}{n}-\nu)}{\Gamma(\frac{x}{n-1}+1)}\mathrm{e}^{-x}
\end{equation}
with equality if and only if $\nu=\frac{1}{n}$.  If $\nu>\frac{1}{n}$, the strict inequality is reversed and holds for all $x>(n-1)(\nu-\frac{1}{n})$.

(ii). Let $x>0$.  Then, for $0\leq \nu\leq\frac{1}{2}$, we have
\begin{equation}\label{label2}K_{\nu}(x)\geq\sqrt{\frac{\pi}{2}}\frac{x^{\nu}\Gamma(x+\frac{1}{2}-\nu)}{\Gamma(x+1)}\mathrm{e}^{-x}
\end{equation}
with equality if and only if $\nu=\frac{1}{2}$.  If $\nu>\frac{1}{2}$, the strict inequality is reversed and holds for all $x>\nu-\frac{1}{2}$.
\end{theorem}

\begin{proof}We first note that part (ii) follows immediately from setting $n=2$ in part (i), due to the relation (\ref{rel}).  In order to establish part (i), we recall the integral representation of the kernel $\lambda_\nu^{(n)}$:  
\[\lambda_\nu^{(n)}(x)=\frac{(2\pi)^{(n-1)/2}\sqrt{n}\big(\frac{x}{n}\big)^{n\nu}}{\Gamma(\nu+1-\frac{1}{n})}\int_1^\infty(t^n-1)^{\nu-\frac{1}{n}}\mathrm{e}^{-xt}\,\mathrm{d}t.\]
Setting $t=\frac{2}{n-1}u+1$ gives
\begin{equation}\label{asdf}\lambda_{\nu}^{(n)}(x)=\frac{(2\pi)^{(n-1)/2}\sqrt{n}\big(\frac{x}{n}\big)^{n\nu}}{\Gamma(\nu+1-\frac{1}{n})}\frac{2}{n-1}\mathrm{e}^{-x}\int_0^{\infty}\bigg(\bigg(\frac{2u}{n-1}+1\bigg)^n-1\bigg)^{\nu-\frac{1}{n}}\mathrm{e}^{-\frac{2x}{n-1}u}\,\mathrm{d}u.
\end{equation}
We now suppose that $0\leq \nu<\frac{1}{n}$ and prove that under this condition inequality (\ref{label}) is strict.  
%For $u>0$ we have $\mathrm{e}^{2u}-1=\sum_{k=1}^{\infty}\frac{(2u)^k}{k!}>2u+2u^2$, and so, for $x>0$,
For $u>0$, we have 
\begin{align*}\frac{n}{n-1}\big(\mathrm{e}^{2u}-1\big)&=\frac{n}{n-1}\sum_{k=1}^{\infty}\frac{(2u)^k}{k!}\\
&>\sum_{k=1}^n\frac{(2u)^k}{k!}\frac{n}{n-1}\times\frac{n-1}{n-1}\times \cdots\times\frac{n-k+1}{n-1} \\
&=\sum_{k=1}^n\frac{(2u)^k}{k!}\frac{n!}{(n-k)!}\cdot\frac{1}{(n-1)^k} \\
&=\sum_{k=1}^n\binom{n}{k}\bigg(\frac{2u}{n-1}\bigg)^k=\bigg(1+\frac{2u}{n-1}\bigg)^n-1.
\end{align*}
Applying this inequality to (\ref{asdf}) yields
\begin{align}\label{first}\lambda_{\nu}^{(n)}(x)&>\frac{(2\pi)^{(n-1)/2}\sqrt{n}\big(\frac{x}{n}\big)^{n\nu}}{\Gamma(\nu+1-\frac{1}{n})}\frac{2}{n-1}\bigg(\frac{n}{n-1}\bigg)^{\nu-\frac{1}{n}}\mathrm{e}^{-x}\int_0^{\infty}\mathrm{e}^{-2xu}(\mathrm{e}^{2u}-1)^{\nu-\frac{1}{n}}\,\mathrm{d}u\\
&=\frac{(2\pi)^{(n-1)/2}\sqrt{n}\big(\frac{x}{n}\big)^{n\nu}}{\Gamma(\nu+1-\frac{1}{n})}\frac{2}{n-1}\bigg(\frac{n}{n-1}\bigg)^{\nu-\frac{1}{n}}\mathrm{e}^{-x}\int_0^{\infty}\mathrm{e}^{-(\frac{2x}{n-1}+\frac{2}{n}-2\nu)u}(1-\mathrm{e}^{-2u})^{\nu-\frac{1}{n}}\,\mathrm{d}u.\nonumber
\end{align}
Making the change of variables $y=\mathrm{e}^{-2u}$ gives
\begin{align*}\int_0^{\infty}\mathrm{e}^{-(\frac{2x}{n-1}+\frac{2}{n}-2\nu)u}(1-\mathrm{e}^{-2u})^{\nu-\frac{1}{n}}\,\mathrm{d}u &=\frac{1}{2}\int_0^1(1-y)^{\nu-\frac{1}{n}}y^{\frac{x}{n-1}+\frac{1}{n}-\nu-1}\,\mathrm{d}y \\
&=\frac{1}{2}B(\nu+1-\tfrac{1}{n},\tfrac{x}{n-1}+\tfrac{1}{n}-\nu) \\
&=\frac{\Gamma(\nu+1-\frac{1}{n})\Gamma(\frac{x}{n-1}+\frac{1}{n}-\nu)}{2\Gamma(\frac{x}{n-1}+1)}, 
\end{align*}
where $B(a,b)$ is the beta function, and we used the standard formula $B(a,b)=\frac{\Gamma(a)\Gamma(b)}{\Gamma(a+b)}$.  This completes the proof that inequality (\ref{label}) holds for $0\leq \nu<\frac{1}{n}$.  When $\nu>\frac{1}{n}$ inequality (\ref{first}) is reversed and so inequality (\ref{label}) is also reversed.  Note that when $\nu>\frac{1}{n}$ the integral in inequality (\ref{first}) only exists if $x>(n-1)(\nu-\frac{1}{n})$. Finally, we note that we have equality when $\nu=\frac{1}{n}$, because (\ref{first}) becomes an equality in this case.
\end{proof}

\begin{corollary}Let $0\leq\nu<\frac{1}{2}$.  Then for all $x>0$,
\[\bigg(\frac{x}{x+\frac{1}{2}-\nu}\bigg)^{\nu+\frac{1}{2}}<\sqrt{\frac{2}{\pi}}\mathrm{e}^xK_{\nu}(x)<1.\]
\end{corollary}

\begin{proof}The upper bound holds because $K_\nu(x)<K_{1/2}(x)=\sqrt{\frac{\pi}{2x}} \mathrm{e}^{-x}$ for all $\nu<\frac{1}{2}$ (see \cite{ifantis}).  The lower bound follows from Theorem \ref{beso} and an application of the inequality $\frac{\Gamma(x+a)}{\Gamma(x+1)}>\frac{1}{(x+a)^{1-a}}$ for $0<a<1$ (see \cite{gautschi}).
\end{proof}

\begin{remark}The following bounds for $K_{\nu}(x)$ were obtained by \cite{luke}:
\[1-\frac{\frac{1}{2}(\frac{1}{4}-\nu^2)}{x+\frac{1}{2}(\frac{1}{4}-\nu^2)}<\sqrt{\frac{2x}{\pi}}\mathrm{e}^xK_{\nu}(x)<1-\frac{\frac{1}{2}(\frac{1}{4}-\nu^2)}{x+\frac{1}{4}(\frac{9}{4}-\nu^2)}, \quad x>0, \:\: 0\leq \nu<\tfrac{1}{2}.\]
Despite taking a relatively simple form, numerical experiments show that, for $0\leq\nu<\frac{1}{2}$, the bounds of \cite{luke} and Theorem \ref{beso}, part (ii) are remarkably accurate for all but very small $x$, for which the modified Bessel function $K_\nu(x)$ has a singularity as $x\downarrow0$.  The bound (\ref{label2}) outperforms the lower bound of \cite{luke} for very small $x$, as it is $O(x^\nu)$ as $x\downarrow0$, as opposed to $O(x^{1/2})$ which is the case for that bound of \cite{luke}.  However, the bound of \cite{luke} performs better for large $x$.     
\end{remark}

\section*{Acknowledgements}

The author is supported by a Dame Kathleen Ollerenshaw Research Fellowship.

\end{document}